\documentclass[11pt]{article}
\usepackage{graphicx}
\usepackage{amssymb}
\usepackage{epstopdf}
\usepackage{amsmath}
\usepackage{amsfonts}
\newcommand{\braket}[2]{\langle #1,#2 \rangle}
\newcommand{\la}{\lambda}

\DeclareSymbolFont{AMSb}{U}{msb}{m}{n}
\DeclareMathSymbol{\N}{\mathbin}{AMSb}{"4E}
\DeclareMathSymbol{\Z}{\mathbin}{AMSb}{"5A}
\DeclareMathSymbol{\R}{\mathbin}{AMSb}{"52}
\DeclareMathSymbol{\Q}{\mathbin}{AMSb}{"51}
\DeclareMathSymbol{\I}{\mathbin}{AMSb}{"49}
\DeclareMathSymbol{\C}{\mathbin}{AMSb}{"43}

\begin{document}

\addtolength{\textheight}{0 cm}
\addtolength{\hoffset}{0 cm}
\addtolength{\textwidth}{0 cm}
\addtolength{\voffset}{0 cm}

\setcounter{secnumdepth}{5}
 \newcommand{\1}{I\hspace{-1.5 mm}I}
\newtheorem{proposition}{Proposition}[section]
\newtheorem{theorem}{Theorem}[section]
\newtheorem{lemma}[theorem]{Lemma}
 \newtheorem{coro}[theorem]{Corollary}
\newtheorem{remark}[theorem]{Remark}
\newtheorem{ex}[theorem]{Example}
\newtheorem{claim}[theorem]{Claim}
\newtheorem{conj}[theorem]{Conjecture}
\newtheorem{definition}[theorem]{Definition}
 
\newtheorem{application}{Application}
 
\newtheorem{corollary}[theorem]{Corollary}

\def\LX{{\cal L}(X)}
\def\LY{{\cal L}(Y)}
\def\LH{{\cal L}(H)}
 \def\ASD{{\cal L}_{\rm AD}(X)}
 \def\ASDY{{\cal L}_{\rm AD}(Y)}
\def\ASDH{{\cal L}_{\rm AD}(H)}
 \def\ASDP{{\cal L}^{+}_{\rm AD}(X)}
  \def\ASDYP{{\cal L}^{+}_{\rm AD}(Y)}
   \def\ASDHP{{\cal L}^{+}_{\rm AD}(H)}
 
\def\CX{{\cal C}(X)}
\def\CY{{\cal C}(Y)}
\def\CH{{\cal C}(H)}
 
\def\PX{{\cal A}(X)}
\def\PY{{\cal A}(Y)}
\def\PH{{\cal A}(H)}
\def\phi{{\varphi}}
\def\AH{A^{2}_{H}}
 \def\H{{\cal H}}

\title{ Iterations of anti-selfdual Lagrangians and applications to Hamiltonian systems and multiparameter gradient flows}
\author{ Nassif  Ghoussoub\thanks{Research partially supported by a grant
from the Natural Sciences and Engineering Research Council of Canada.} \quad and \quad Leo Tzou
 \\
\small Department of Mathematics,
\small University of British Columbia, \\
\small Vancouver BC Canada V6T 1Z2 \\
\small {\tt nassif@math.ubc.ca} \\
\small {\tt leo@pims.math.ca}
\\
%\today\\ 
%{January 20, 2004}\\
}

\maketitle

\section*{Abstract} Anti-selfdual Lagrangians on a state space lift to path space provided one adds a suitable selfdual boundary Lagrangian.   This process can be iterated by considering the path space as a new state space for the newly obtained anti-selfdual Lagrangian. We give here two applications for these remarkable permanence properties. In the first, we  establish for certain convex-concave Hamiltonians ${\cal H}$ on a --possibly infinite dimensional--symplectic space $H^2$,  the existence of a solution for the Hamiltonian system $-J\dot u (t)=\partial {\cal H} (u(t))$ that connects in a given time $T>0$, two Lagrangian submanifolds. Another application deals with the construction of a multiparameter  gradient flow for a convex potential.  Our methods are based on the new  variational calculus for anti-selfdual Lagrangians developed in \cite{G2}, \cite{G3} and \cite{GT1}.  
 
 \section{Introduction} Given two convex and lower semi-continuous functions $(\phi_1, \phi_2)$ on $\R^n$, we consider the  Hamiltonian  $\H$ on $\R^{2n}$ defined by $\H(x,y) = \phi_1(x) - \phi_2(y)$ and we look for solutions for the Hamiltonian system $-J\dot u (t)=\partial {\cal H} (u(t))$ that connects in time $T>0$, the Lagrangian submanifolds 
\[
\hbox{ $ L_1=\{(x, y)\in \R^{2n}; -y \in A_1x+\partial \psi_1(x)\}$\quad  to \quad $L_2=\{(x, y)\in \R^{2n}; y\in A_2x+\partial \psi_2(x)\}$. }
\]
where $\psi_1, \psi_2$ are convex lower semi-continuous functions on $\R^n$ and $A_1$, $A_2$ are positive (but not neccesarily self-adjoint) matrices.  
 In other words, we are looking for a solution on $[0,T]$ for the  Hamiltonian system:
 
 \begin{eqnarray}
{\dot x} (t) &\in& \partial_2\H(x(t),y(t)) \nonumber\\
-{\dot y}(t) &\in& \partial_1\H(x(t),y(t))
 \end{eqnarray}
 
 with the following boundary conditions
 \begin{equation}
-y(0)- A_1x(0)\in \partial  \psi_1 (x(0))\quad {\rm and}\quad 
   y(T)- A_2x(T)\in \partial  \psi_1 (x(T)).
\end{equation}
We shall show that a solution can indeed be obtained  by minimizing the following functional
 {  \small
  \[I(x,y) =   \int_0^T \Phi((x(t),y(t)) + \Phi^*(-\dot{y}(t),-\dot{x}(t))dt + \psi_1( x(0))+ \psi_1^*(-y(0)-A_1x(0))+ \psi_2(x(T)) +  \psi_2^*(y(T)-A_2x(T)).
\]
  }
on the space $A^2([0,T]; \R^{2n}) = \{ u=(x,y):[0,T] \rightarrow \R^{2n};   \dot{u} \in L^{2}_{\R^{2n}}  \}$,
where here $\Phi$ is the convex function $\Phi (x,y)=\phi_1(x) + \phi_2(y)$ on $\R^{2n}$ and $\Phi^*$ is its Legendre transform. The equation is obtained from the fact that the infimum is actually $0$, which is the main point of the exercise.  

Actually, this is a particular case of a much more general result. For one, the method is infinite dimensional and $\R^n$ can be replaced by any Hilbert space $H$ and for PDE puposes, the domain can be an evolution pair $X\subset  H \subset X^*$ where $X$ is a Banach space dense in $H$.  More importantly, the theorem is really about the existence of a path connecting in prescribed time $T$,  two given  ``anti-selfdual" Lagrangian submanifolds in $H^2$ through an "anti-selfdual" Lagrangian submanifold in phase space $H^4$. 
Let us first recall the following notions from \cite{G2}. 

\begin{definition} (1) A convex lower semi-continuous functional $L: H\times H \to \R \cup\{+\infty\}$ (resp., $\ell: H\times H \to \R \cup\{+\infty\})$ is said to be $R$-antiselfdual (resp., $R$-selfdual) for some automorphism $R:H\to H$ if 
\[
\hbox{ $L^*(p,x)=L(-Rx, -Rp)$ \quad (resp., $\ell^*(x,p)=\ell (-Rx, Rp)$) for any $(x,p)\in H\times H$.}
\]
(2) An $R$-antiselfdual manifold $M$ in $H\times H$ is a set of the form 
\[
M=\{(x,p)\in H\times H;\,\, L(x,p)+\langle Rx, p\rangle =0\}
\]
where $L$ is an $R$-antiselfdual Lagrangian on $H$. 
\end{definition}
Typical examples are 
\[
M_{+, \psi}=\{(x,p)\in H\times H; \psi (x)+\psi^*(-p)+\langle x, p\rangle =0\}=\{(x,p)\in H\times H; \, p\in -\partial \psi (x)\}. 
\]
and 
\[
M_{-, \psi}=\{(x,p)\in H\times H; \psi (x)+\psi^*(p)-\langle x, p\rangle =0\}=\{(x,p)\in H\times H; \, p\in \partial \psi (x)\}. 
\]
where $\psi$ is a convex lower semi-continous function on $H$ and where $R(x)=x$ for $M_{+, \psi}$ and $R(x)=-x$ for $M_{-, \psi}$. 

Moreover, if $A:H\to H$ is a bounded skew-adjoint operator on $H$, then the following manifolds are also $(+I)-ASD$ (resp., $(-I)-ASD)$ (See \cite{G2}).
\[
M_{+, \psi, A}=\{(x,p)\in H\times H; \psi (x)+\psi^*(-Ax-p)+\langle x, p\rangle =0\}=\{(x,p)\in H\times H;\,  -p\in (A+\partial \psi) (x)\}. 
\]
and 
\[
M_{-, \psi, A}=\{(x,p)\in H\times H; \psi (x)+\psi^*(-Ax+p)-\langle x, p\rangle =0\}=\{(x,p)\in H\times H;\,  p\in (A+\partial \psi) (x)\}. 
\]
 The condition that $A$ is a skew-adjoint operator can be replaced by the hypothesis that  it is merely positive, i.e., that $\langle A x,x\rangle \geq 0$ for every $x\in H$. Indeed, one can decompose $A$ into its symmetric part $A^s=\frac{1}{2}(A x+A^*x)$ and its skew-symmetric part  $A^a=\frac{1}{2}(A x-A^*x)$. Then, the manifold 
 \[
 {M}_{+,\psi, A} = \{(x, p)\in H\times H;  \, -p - Ax \in\partial\psi(x)\}
 \]
 is equal to the $(+I)$-ASD manifold 
  \[
 {\cal M}_{+,\tilde \psi, A^a} = \{(x,p)\in H\times H;  \, -p - A^a x \in\partial \tilde \psi(x)\}
 \]
where $\tilde \psi(x) =\psi(x) +\frac{1}{2}\langle Ax,x\rangle$, while the manifold 
 \[
 {M}_{-,\psi, A} = \{(x, p)\in H\times H;  \, p - Ax \in\partial\psi(x)\}
 \]
 is equal to the $(-I)$-ASD manifold 
  \[
 {\cal M}_{-,\tilde \psi, A^a} = \{(x,p)\in H\times H;  \, p - A^a x \in\partial \tilde \psi(x)\}
 \]
This will allow us --in the sequel-- to reduce many of the proofs for statements concerning bounded positive operators  to the case where they are skew adjoint.

Consider now a convex lower semi-continuous function $\Phi$ on $H\times H$ and let $S: H\times H \to H\times H$ be the automorphism $S(p,q)=(q,p)$, then one can easily check that the following manifold 
\[
{\cal M}_{S, \Phi}:= \{\big((x_1, x_2),(p_1,p_2)\big)\in H^2\times H^2;\,  \big(-p_2, -p_1\big)\in\big(\partial_1\Phi(x_1,x_2),\partial_2\Phi(x_1,x_2)\big)\}
\]
is $S$-antiselfdual, and can be written as
\[
{\cal M}_{S, \Phi}:= \{\big((x_1, x_2),(p_1,p_2)\big)\in H^2\times H^2; \,  \Phi(x_1,x_2)+\Phi^*(-S(p_1,p_2))+ \langle(x_1,x_2), S(p_1,p_2)\rangle=0\}
\]

Our main theorem in section 2 below asserts that under very general conditions, one should be able for any time $T>0$,  to connect any given  $(+I)$-ASD submanifold in $H^2$ to a given  $(-I)$-ASD submanifold in $H^2$ through a path in phase space $(x(t), \dot x(t)) $ that lies on a given $S$-ASD submanifold in $H^4$.

The proof relies on the extremely useful fact that if $L$ is an $R$-antiselfdual Lagrangian on state space and if $\ell$ is an $R$-selfdual boundary Lagrangian then the following Lagrangian defined by 
\begin{equation}
\label{iteration}
{\cal L}(x,p) := 
\begin{cases}
{\int_0^T L(t,x(t), \dot x (t)+ p(t))dt +  \ell (x(0), x(T)) \mbox{ if}\, \,  \dot x\in L^2_H}\cr{+\infty \quad \quad \quad \quad \quad \quad \quad \quad \quad \quad \quad \quad \mbox{ elsewhere}}\cr
\end{cases}
\end{equation}
 is also an $R$-antiselfdual Lagrangian on path space $L^2_H[0,T]$.

 In section 3, we exploit the antiselfduality of this new Lagrangian to lift it  to another ASD Lagrangian on  a new path space $L^2([0,S]; L^2_H([0,T]$. Applied to  the basic ASD Lagrangian $L(x,p)=\phi (x) +\phi^*(-p)$  associated to a given convex lower semi-continuous function $\phi$,  this leads to the construction  for any $x_0\in H$, $T>0$ and $S>0$, of surfaces $\hat x(t,s)$ verifying 
 for almost all  $(s,t)\in[0,S]\times[0,T]$
  \[
\frac{\partial \hat x}{\partial t}(s,t) +\frac{\partial \hat x}{\partial s}(s,t)\in -\partial \phi (\hat x(s, t))
\]
\[\hat x(0,t) = x_0\ {\rm a.e.} \ t\in [0,T]\]
\[\hat x(s,0) = x_0\ {\rm a.e.} \ s\in [0,S].
\]
It is clear that this process can be iterated to obtain some kind of a multiparameter gradient flow for any convex potential.

\section{Connecting Lagrangian submanifolds}
As mentioned above, the key ingredient in what follows is the fact that if $L$ is an $R$-ASD Lagrangian on a space $H$, then  --under suitable boundedness conditions-- the Lagrangian $\cal L$ defined in (\ref{iteration}) is then R-ASD on the path space $L^2_H$. The proof of the main result in this section  requires however that $\cal L$ be only {\it partially $R$-antiselfdual} on path space (See \cite{G2}) which holds --as proved below-- without additional boundedness conditions. The infinite dimensional framework required by the applications to PDE can be formulated in many settings. We describe some of them in varying levels of detail. 

\subsection{The Hilbertian framework}

Let $H$ be a Hilbert space with $\braket \,{\,}\,{\,}$ as scalar product and let $[0,T]$ be a fixed real interval.   For $\alpha \in (1, +\infty)$, we consider the classical space  
$L^{\alpha}_{H}$ of Bochner integrable functions from $[0,T]$ into $H$ with norm denoted by $\|\cdot\|_\alpha$, as well as the reflexive Banach space
 $
A^{\alpha}_{H} = \{ u:[0,T] \rightarrow H; \,\dot{u} \in L^{\alpha}_{H}  \}
$
consisting of all absolutely continuous arcs $u: [0,T]\to H$, 
equipped with the norm
\[
     \|u\|_{A^{^{\alpha}}_{H}} = \|u(0)\|_{H} +
    ( \int_0^T \|\dot{u}\|^{\alpha} dt)^{\frac{1}{\alpha}}.
\]
It is clear that $A^{^{\alpha}}_{H}$ can be identified with the product space
$H
\times L^{\alpha}_{H}$, and that its dual  $(A^{\alpha}_{H})^*$ can also be
identified with $ H \times
L^{\beta}_{H}$ (where $\frac{1}{\alpha} +\frac{1}{\beta}=1$)  via the formula:
\[
     \braket{u}{(a,p)}_{_{A^{\alpha}_{H},H \times L_{H}^{\beta}}} =
      \braket{u(0)}{a}_{_H} + \int_0^T \braket{\dot{u}(t)}{p(t)}dt.
\]
We consider the following action functional on  $A^{\alpha}_{H}$:
\[
I_{\ell,L}(u) = \int_0^T L(t, u(t),\dot{u}(t))dt + \ell (u(0),u(T))
\]
where
\[ \ell: H
\times
H
\rightarrow
{\R \cup \{+\infty\}}\quad \hbox{\rm and}\quad 
     L:[0,T]\times H\times H
        \rightarrow {\R \cup \{+\infty\}}
\]
are two appropriate Lagrangians. We shall always assume that $L$ is measurable with respect to   the  $\sigma$-field in $[0,T]\times H\times H$ generated by the products of Lebesgue sets in  $[0,T]$ and Borel sets in $H\times H$, and  that  $\ell$ and $L(t,\cdot,\cdot)$ are convex, lower semi-continuous valued in $\R \cup \{+\infty\}$ but not identically $+ \infty$. 

\begin{theorem} \label{existence}
Assume that $R$ is an automorphism of $H$, that $L (t, \cdot, \cdot): H\times H \to \R \cup \{+\infty\}$ is $R$-antiselfdual for each $t\in [0,T]$ and that  $\ell$ is $R$-selfdual.  Assume
\begin{equation}
\label{boundedness}
\hbox{ $L(t,y,0)\leq C(1+\Vert y\Vert_H^\beta)$ for $y\in H$  and $a\mapsto \ell(a,0)$  is bounded on the bounded sets of $H$.}
\end{equation}
 Then there exists $\hat x \in A^\alpha_H$ such that
\begin{equation}
\label{min=zero}
I_{L,\ell}(\hat x ) = \inf\limits_{x \in A^\alpha_H} I_{L,\ell}(x) = 0.
\end{equation}
 \end{theorem}

For the proof, we consider the  functional $J^\alpha_{L,\ell}: (A^\alpha_H)^* \cong H\times L^\beta_H \to\R \cup \{+\infty\}$ defined by:
\[
J^\alpha_{L,\ell}(a,y(\cdot)) := \inf\limits_{x(\cdot)\in A^\alpha_H}\{ \int_0^T L(t,x(t)+ y(t),\dot x(t))dt +\ell(x(0)+a,x(T))\}.
\]
The key to the proof is the following proposition
\begin{proposition} Assume that $R$ is an automorphism of $H$, that $L (t, \cdot, \cdot): H\times H \to\R \cup \{+\infty\}$ is $R$-antiselfdual for each $t\in [0,T]$ and that  $\ell$ is $R$-selfdual. Then 
\begin{enumerate}
\item The functional $J^\alpha_{L,\ell}$ is convex on $H\times L^\beta_H$ and its Legendre transform  in the duality
$(H\times L^{\beta}_{H}, A_{H}^{\alpha})$ satisfies  for any  $x\in A^\alpha_H$, 
 \begin{equation}
 \label{invariance}
(J^\alpha_{L,\ell})^*(x) = \int_0^T L(t,-Rx(t),-R\dot x(t))dt +\ell(-Rx(0),-Rx(T))=I_{\ell,L}(-Rx). 
\end{equation}
\item If $J^\alpha_{L,\ell}$ is subdifferentiable at $(0,0)$ on the space $H\times L^\beta_H$, then there exists $\hat x \in A^\alpha_H$ such that
$
I_{L,\ell}(\hat x ) = \inf\limits_{x \in A^\alpha_H} I_{L,\ell}(x) = 0.$
\end{enumerate}
 \end{proposition}
  \noindent {\bf Proof:}  1) The convexity of  $J_{\ell,L}$ is easy to establish.    Fix now $p \in A^{\alpha}_{H}$ and write:
\[
     J^*_{\ell, L}(p) = \sup_{a \in H} \sup_{y \in L_H^{\beta}} \sup_{u \in
     A^{\alpha}_{H}} \left\{ \langle a,p(0) \rangle + \int_0^T \left[ \braket{y(t)}{\dot{p}(t)}
     - L(t, u(t)+y(t),\dot{u})\right] dt -\ell (u(0)+a,u(T)) \right\}.
\]
 Make a substitution $
u(0) + a = a' \in H \quad \hbox{\rm and $u+y=y' \in L_H^\beta$}$, 
we obtain
\[
     J^*_{\ell,L} (p) = \sup_{a' \in H} \sup_{y' \in L_H^{\beta}}
     \sup_{u \in A_{H}^{^\alpha}} \left\{ \langle a'-u(0),p(0) \rangle - \ell (a',u(T)) +
     \int_0^T\left[ \braket{y'(t)-u(t)}{\dot{p}(t)} - L(t,y'(t),\dot{u}(t)\right] dt  \right\}.
\]
Since $\dot u \in L_H^\alpha$ and $u\in L_H^\beta$, we have $
\int_0^T \braket{u}{\dot{p}}=- \int_0^T \braket{\dot u}{p} +
\braket{p(T)}{u(T)} - \braket{p(0)}{u(0)},$
    which implies
\begin{eqnarray*}
     J^*_{\ell ,L} (p) = \sup_{a' \in H} \sup_{y' \in L_H^{\beta}} \sup_{u
\in A_{H}^{^\alpha}}
\{ \langle a',p(0) \rangle &+&
     \int_0^T
\{\braket{y'}{\dot{p}}+\braket{\dot{u}}{p}-L(t,y'(t),\dot{u}(t))\}
dt\\
&-&\braket{u(T)}{p(T)}-\ell \big( a',u(T)\big) \}.
\end{eqnarray*}
It is now convenient to identify $A_{H}^{^\alpha}$ with $H\times L_H^\alpha$
via the correspondence:
 $    (c,v) \in H\times L_H^\alpha \mapsto  c+\int_t^T v(s)\, ds\in
A_{H}^{^\alpha}$ and $
     u\in A_{H}^{^\alpha}   \mapsto  \big( u(T),-\dot u(t)\big)\in H\times
L_H^\alpha$. 
We finally obtain
\begin{eqnarray*}
     J^*_{\ell ,L}(p)&=& \sup_{a' \in H} \sup_{c\in H}
\left\{ \langle a',p(0) \rangle +  \braket{-c}{p(T)}-\ell(a',c)\right\}\\
&&+\sup_{y'\in L_H^\beta}\sup_{v\in L_H^\alpha}\left\{
\int_0^T\left[ \braket{y'}{\dot{p}} +\braket{v}{p} - L(t,y'(t),v(t))\right] dt\right\}\\
   &=& \int_0^T L^*(t, \dot p (t), p(t)) dt +\ell^* (p(0), -p(T))\\
   &=& \int_0^T L(t, -R p (t), -R\dot p(t)) dt +\ell (-Rp(0), -Rp(T))\\
  &=& I_{\ell, L}(-Rp).
\end{eqnarray*}
 
 2) Since $R$ is an automorphism, weak duality gives
\[
\inf_{u\in A_{H}^2} I_{\ell ,L}(u)\geq 
 \sup_{A_{H}^2}  -
J_{\ell ,L}^* (u)= \sup_{A_{H}^2}  -
I_{\ell ,L} (-Ru)=\sup_{A_{H}^2}  -
I_{\ell ,L} (u)=-\inf_{u\in A_{H}^2} I_{\ell ,L}(u)
\]
and $\inf_{u\in A_{H}^2} I_{\ell ,L}(u)$ is therefore non negative. 

On the other hand, if we pick
 $\hat x \in \partial J^\alpha_{L,\ell} (0,0)$, we get
 \[
-  \inf_{A_{H}^2} I_{\ell ,L}(u)=-J^\alpha_{L,\ell}(0)=(J^\alpha_{L,\ell})^*(\hat x)=I_{\ell ,L}(-R \hat x)\geq   \inf_{A_{H}^2} I_{\ell ,L}(u)
\]
which means that $ \inf_{A_{H}^2} I_{\ell ,L}(u) \leq 0$. It follows that 
$
\inf_{A_{H}^2} I_{\ell ,L}(u)=I_{\ell ,L}(-R \hat x)=0.
$\\

 \noindent{\bf Proof of Theorem 2.1}:  It remains to show  that  the convex functional $J_{\ell ,L}$ is sub-differentiable  at $(0,0)$ on the space $H\times L^\beta_H$ so as to conclude using  Proposition 2.1. But the boundedness assumptions  (\ref{boundedness}) on $L$ and $\ell$ immediately give
\[
J_{\ell,L}(a,y)
   \leq  \int_0^T L(t, y(t), 0)dt + \ell (a, 0)
 \leq \int_{0}^{T} C (1+\|y(t)\|_{H}^{\beta})dt + \ell (a, 0)
 \]
which means that $J_{\ell,L}$ is bounded on the bounded sets of $H\times L^\beta_H$ and since it is convex, it is therefore subdifferentiable at $(0,0)$.

\begin{theorem}
\label{lagrangian manifold} Let $\psi_1$ and $\psi_2$ be two convex and lower semi-continuous functions on 
 a Hilbert space $E$, let $A_1, A_2 :E\to E$ be  bounded  positive operators and consider the manifolds
 \[
{\cal M}_1:= {\cal M}_{+,\psi_1, A_1} = \{(x_1, x_2)\in E\times E;  \, -x_2- A_1 x_1 \in\partial\psi_1(x_1)\}
 \]
 and
\[
{\cal M}_2:={\cal M}_{-,\psi_2, A_2}= \{(x_1, x_2)\in E\times E; \,  x_2 - A_2 x_1\in\partial\psi_2(x_1)\}.
\]
Let $\Phi: [0,T]\times K\times K\to \R$ be such that $\Phi (t, \cdot, \cdot) $ is convex and lower semi-continuous  for  each $t\in [0,T]$ and consider  the evolving manifold
\[
{\cal M}_3(t):={\cal M}_{S,\Phi}(t) = \{\big((x_1, x_2),(p_1,p_2)\big)\in E^2\times E^2; \, (-p_2, -p_1)\in\big(\partial_1\Phi(t,x_1,x_2),\partial_2\Phi(t,x_1,x_2)\big)\}
\]
Now assume that $\psi_1$ is coercive and bounded on bounded sets of $E$, $\psi_2$ is bounded below with $0$ in its domain, while for every $ t \in[0,T]$ we have
 \begin{equation}
\label{interior.boundedness}
\hbox{ 
$\Phi(t,x_1,x_2)\leq C(1+\Vert x_1\Vert_E^\beta+\Vert x_2\Vert_E^\beta).$}
\end{equation}
 Then there exists $x \in A^\alpha_{E\times E}$ such that:
 \[
 \hbox{ $x(0)\in {\cal M}_1$,\, $x(T)\in {\cal M}_2$\, and\,  $(x(t), \dot{x}(t)) \in {\cal M}_3(t)$ for a.e. $t\in [0,T]$.}
 \]
\end{theorem}
We shall need the following easy but interesting lemma.
\begin{lemma} Suppose $\ell_1$ (resp., $\ell_2$) is an (+I)-anti-selfdual  Lagrangian (resp., an (+I)-anti-selfdual  Lagrangian on the Hilbert space $E\times E$, then  the Lagrangian  $\ell:E^2\times E^2\to \R$ defined by 
\[
\ell ((a_1,a_2), (b_1, b_2))=\ell_1(a_1,a_2)+\ell_2 (b_1, b_2)
\]
 is $S$-selfdual on $E^2\times E^2$ where $S$ is the automorphism on $E\times E$ defined by $S(x_1, x_2) = (x_2, x_1)$.  
 
 In particular, if  $\psi_1$ and $\psi_2$ are convex lower semi-continuous on $E$ and if $A_1,\  A_2$ are bounded skew-adjoint operators on $E$,  then the Lagrangian $\ell(\cdot,\cdot): H\times H$ defined by 
\[\ell(a,b) := \psi_1(a_1) + \psi_1^*(-A_1a_1-a_2) + \psi_2(b_1) + \psi_2^*(-A_2b_1 +b_2)\]
 is S-selfdual.  
\end{lemma}
The proof is left to the interested reader (See also \cite{G2}).\\

\noindent{\bf Proof of Theorem 2.2:} In view of Remark 1.2, we can assume without loss that $A_1$ and $A_2$ are skew adjoint operators. Let $H=E\times E$ and 
consider the S-anti-selfdual Lagrangian on $H\times H$ defined by
\[L(t,x,p) := \Phi(t,x) + \Phi^*(t,-Sx)\]
 as well as the $S$-selfdual boundary Lagrangian $\ell: H\times H$ defined by 
\[
\ell(a,b) := \psi_1(a_1) + \psi_1^*(-A_1a_1-a_2) + \psi_2(b_1) + \psi_2^*(-A_2b_1 +b_2).
\]
Since  $\Phi(t,x_1,x_2)\leq C(1+\Vert x_1\Vert_K^\beta+\Vert x_2\Vert_K^\beta)\forall t \in[0,T]$ and since $\phi$ is coercive and bounded on bounded sets,  the functional defined on $A^\alpha_{E\times E}[0,T]$ by
{\small
\[
I_{L,\ell}(u)=\int_0^T (\Phi(t,u(t)) + \Phi^*(t,-S\dot{u}(t)))dt
 + \psi_1({u_1}(0)) + \psi_1^*(-A_1{u_1}(0)-{u_2}(0)) + \psi_2({u_1}(T)) + \psi_2^*(-A_2{u_1}(T) +{u_2}(T))
\]
}
satisfies all the hypothesis of  Theorem \ref{existence}. Hence there exists $x(\cdot)\in A^\alpha_{E\times E}$ such that $I_{L,\ell}(x)= 0$. Therefore,
\begin{eqnarray*}
0&=&\int_0^T \Phi(t,x(t)) + \Phi^*(t,-S\dot{x}(t))dt\\
&& + \psi_1({x_1}(0)) + \psi_1^*(-A_1{x_1}(0)-{x_2}(0)) + \psi_2({x_1}(T)) + \psi_2^*(-A_2{x_1}(T) +{x_2}(T))\\
&\geq& \int_0^T \braket{x(t)}{-S\dot{x}(t)}dt + \psi_1({x_1}(0)) + \psi_1^*(-A_1{x_1}(0)-{x_2}(0)) + \psi_2({x_1}(T)) + \psi_2^*(-A_2{x_1}(T) +{x_2}(T))\\
&=& -\int_0^T \frac{d\braket{x_1(t)}{x_2(t)}}{dt}dt  + \psi_1({x_1}(0)) + \psi_1^*(-A_1{x_1}(0)-{x_2}(0)) + \psi_2({x_1}(T)) + \psi_2^*(-A_2{x_1}(T) +{x_2}(T))\\
&=&\braket{{x_1}(0)}{{x_2}(0)} +\psi_1({x_1}(0)) + \psi_1^*(-A_1{x_1}(0)-{x_2}(0)) \\
&&- \braket{{x_1}(T)}{{x_2}(T)}+ \psi_2({x_1}(T)) + \psi_2^*(-A_2{x_1}(T) +{x_2}(T))\\
&\geq&0.
\end{eqnarray*}
This means that every inequality in this chain is an equality, hence three applications of the limiting case in Legendre-Fenchel duality gives:
\[-\big(\dot x_2(t), \dot x_1(t)\big)\in\big(\partial_1\Phi({x_1}(t),{x_2}(t)),\partial_2\Phi({x_1}(t),{x_2}(t))\big)\, {\rm a.e.}\ t\in[0,T]\]
\[-A_1x_1(0) - x_2(0)\in\partial\psi_1(x_1(0))\]
\[-A_2x_1(T) + x_2(T)\in\partial\psi_2(x_1(T)).\]
In other words, 
$x(\cdot)\in A^2_{E\times E}$ is such that $x(0)\in {\cal M}_1$, $x(T)\in {\cal M}_1$ and $-(x(t), \dot{x}(t)) \in {\cal M}_3(t)$ for a.e. $t\in [0,T]$
\begin{corollary}
\label{hamiltonian system}
Let $E$ be a Hilbert space and ${\cal H}(\cdot,\cdot) : E\times E\to\R$ be a Hamiltonian of the form ${\cal H}(x_1,x_2) = \phi_1(x_1) - \phi_2(x_2)$ where $\phi_1, \phi_2$ are convex lower semi-continuous functions satisfying 
\[
\phi_1(x_1) + \phi_2(x_2) \leq C\big(1 + \Vert x_1\Vert_K^\beta+\Vert x_2\Vert_K^\beta\big).
\]
 Furthermore, let $\psi_1$, $\psi_2$, $A_1$, and $A_2$ be as in Theorem \ref{lagrangian manifold}. Then there exists $({x_1}, {x_2})\in A^\alpha_{E\times E}([0,T])$ such that for almost all $t\in[0,T]$,
\[-\dot x_2(t)\in\partial_1{\cal H}({x_1}(t), {x_2}(t))\]
\[\dot x_1 (t)\in\partial_2{\cal H}({x_1}(t), {x_2}(t))\]
and satisfying the boundary conditions
\[-A_1x_1(0) - x_2(0)\in\partial\psi_1(x_1(0))\]
\[-A_2x_1(T) + x_2(T)\in\partial\psi_2(x_1(T))\]
\end{corollary}
\noindent{\bf Proof:} This is a restatement of Theorem \ref{lagrangian manifold} for $\Phi(x_1,x_2)=\phi_1(x_1) + \phi_2(x_2)$. 
\begin{corollary}
\label{second order}
Let $E$ be a Hilbert space and let $\phi$ be a convex lower semi-continuous function on $E$ satisfying 
$
\phi (x) \leq C\big(1 + \Vert x\|_H^\beta)$. 
Let $\psi_1$, $\psi_2$, $A_1$, and $A_2$ be as in Theorem \ref{lagrangian manifold}. Then there exists $x \in A^\alpha_{E}([0,T])$ such that for almost all $t\in[0,T]$,
\begin{eqnarray*}
\ddot x(t)&\in&\partial \phi (x(t))\\
 - \dot x(0)&\in&\partial\psi_1(x(0))+A_1x(0) \\
 \dot x(T)&\in&\partial\psi_2(x(T))+A_2x(T).
\end{eqnarray*}
\end{corollary}
\noindent{\bf Proof:}  It si enough to apply the above to $\phi_2=\phi$ and $\phi_1(x_1)=\frac{1}{2}\|x_1\|_H^2$. 

 \subsection{The non-Hilbertian case}

In the  infinite dimensional setting --more suitable for applications to PDEs-- we need the framework of an evolution triple $X \subset H \subset X^*$, where   $H$ is a Hilbert space with $\braket{}{}$ as scalar product, and $X$ is a dense vector subspace of $H$, that  is a reflexive Banach space once equipped
with its norm $\| \cdot \|$.   Assuming the canonical injection $X \rightarrow H$, 
continuous, we identify the Hilbert space $H$ with its dual $H^*$ and we
``inject''
$H$ in $X^*$   in such a way that $
     \braket{h}{u}_{X^*,X} = \braket{h}{u}_H$ for all $h \in
H$
and all $u \in X$. 
This injection  is continuous, one-to-one, and $H$ is also dense in $X^*$. 
 In other words,  the dual $X^*$ of $X$ is represented as the completion of $H$ for the dual norm
$\|h\|=\sup\{{\braket{h}{u}}_H; \|u\|_X \leq 1\}$.\\
 We shall consider here evolution equations with two types of initial conditions. The first ones are those  involving bounded operators in the initial conditions, or boundary Lagrangians on the ambiant Hilbert space $H$ such as Hamiltonian systems of the form:
\begin{eqnarray*}
\dot p (t) &\in& \partial_2\H(p(t),q(t))\\
-\dot q(t)& \in& \partial_1\H(p(t),q(t))\\
p(0) = -q(0) &\&&  p(T) = q(T).
\end{eqnarray*}
We would also like to consider more complex initial conditions:
\[\dot p (t) \in \partial_2\H(p(t),q(t))\]
\[-\dot q(t)\in \partial_1\H(p(t),q(t))\]
\[-A_1p(0) - q(0)\in\partial\psi_1(p(0))\]
\[-A_2p(T) + q(T)\in\partial\psi_2(p(T))\]
where $\psi_1$, $\psi_2$ may only be finite on the space $X$.\\
For the first system the spaces to consider are
\[
A^\alpha_{H,X^*} = \{ u:[0,T] \rightarrow X^*; 
 u(0)\in H, \dot{u} \in L^{\alpha}_{X^*}  \}
\]
equiped with the norm $\|u\|_{A^\alpha_{H,X^*}} = \big(\int_0^T \|\dot u(t)\|_{X^*}^\alpha dt\big)^{1/\alpha} + \|u(0)\|_H$ for $1 < \alpha < \infty$.

For the second system we will need the space
\[
A^\alpha_{X^*} = \{ u:[0,T] \rightarrow X^*; 
 \quad u \, \& \, \dot{u} \in L^{\alpha}_{X^*}  \}
\]
equipped with the norm
$
     \|u\|_{A^{\alpha}_{X^*}} = \|u(0)\|_{X^*} +
     (\int_0^T \|\dot{u}\|_{_{X^{*}}}^{\alpha} dt)^{\frac{1}{\alpha}}.
$

Since the proof of existence for both equations is similar in spirit,  we will only show the detailed proof for the second initial value problem. The other case is left to the interested reader.\\
It is clear that $A^\alpha_{X^*}$ is a reflexive Banach spaces that can be
identified with the product space
$X^*
\times L^{\alpha}_{X^*}$, while its dual  $(A^{\alpha}_{X^*})^*
\simeq X\times
L^{\beta}_{X}$ where $\frac{1}{\alpha} + \frac{1}{\beta} = 1$. The duality is
then given by the formula:
\[
   \braket{u}{(a,p)}_{A^{\alpha}_{X^*},X \times L_{X}^{\beta}} =
 \braket{u(0)}{a} + \int_0^T \braket{\dot{u}(t)}{p(t)}dt
\]
where $\braket{\cdot}{\cdot}$ is the duality on $X$, $X^*$ and
$(\cdot,\cdot)$ is the inner product on $H$.

Let $\ell: X^*
\times
X^*
\rightarrow
{\R \cup \{+\infty\}}$ be convex and weak$^*$-lower semi-continuous on
$X^*\times X^*$, and let
$
     L:[0,T]\times X^*\times X^*
        \rightarrow \R \cup \{+\infty\}
$
be measurable with respect to   the  $\sigma$-field in $[0,T]\times X^{*}\times
X^{*}$ generated by the products of Lebesgue sets in $[0,T]$ and Borel sets in
$X^{*}\times X^{*}$, in such a way that for each $t\in [0,T]$, $L(t,
\cdot,\cdot)$ is convex and weak$^*$-lower semi-continuous on $X^*\times X^*$.
\begin{definition}
Let $R : X^* \to X^*$ be any map. We say that  $L$ is $R$-anti-self-dual and $\ell$ is $R$-selfdual on $X$
if  for all $(p,s)\in X^*\times X^*$, we have
 \[
(\ell |_{X\times
X})^*(p,s)=\ell (-Rp,Rs)  \quad \hbox{\rm and $
(L_{t}|_{X\times X})^*(t, p,s)=L(t, -Rs,-Rp)$}.
\]
where $(L_{t}|_{X\times X})^*$ and $(\ell|_{X\times X})^*$ denote the
Legendre duals of the restrictions of $L_{t}=L(t, \cdot, \cdot)$ and $\ell$ to
$X\times X$. 
\end{definition}

To any such a pair, we
associate the action functional on  $A^{\alpha}_{X^*}$ by:
\[
I_{\ell,L}(u) = \int_0^T L(t, u(t),\dot{u}(t))dt + \ell (u(0),u(T))\]
as well as the corresponding ``variation function'' $J^\alpha_{\ell,L}$
defined on
$(A^{\alpha}_{X^*})^* = X \times L_X^{\beta}$ by
\[
J^\alpha_{\ell,L} (a,y) = \inf \{ \int_0^T L(t, u+y,\dot{u})dt + \ell
(u(0)+a,u(T))\ ;
\ u \in A^{\alpha}_{X^*} \}
\]
 
\begin{theorem} Suppose that $R : X^* \to X^*$ is an automorphism whose restriction to $H$ and 
\label{theorem2}
$X$ is also an automorphism on these spaces. Suppose that for each $t\in [0,T]$, the Lagrangians $L(t,
\cdot)$ and $l$ are two proper convex and
weak$^*$-lower
semi-continuous functions on $X^*\times X^*$ such that $L$ is $R$-anti-self-dual and $\ell$ is $R$-selfdual on $X$. Suppose that for some  $\alpha \in (1, 2]$, $J^\alpha_{\ell ,L}:
X\times L_X^{\beta} \rightarrow {\R \cup \{+\infty\}}$ is
sub-differentiable at
$(0,0)$, then there
exists
$v \in A_{X^*}^\alpha$ such that:
$
 (v(t),\dot v(t))\in  {\rm Dom} (L)$ for almost all
$t\in [0,T]$, and
$
   I_{\ell ,L}(v)=\inf\limits_{A_{X^*}^\alpha}I_{\ell, L}(u)=0.
$
\end{theorem}
Theorem \ref{theorem2} can be proved just like Theorem 2.1 above. The only serious change occurs in the following lemma 
whose proof we include.  

\begin{lemma} Under the above conditions, we have
\label{J*}
$ J^*_{\ell,L}(p) \geq I_{\ell, L} (-Rp)$ for all $p \in
A^{\alpha}_{X^*}$.
\end{lemma}
{\bf Proof:} For $p \in A^{\alpha}_{X^*}$, write:
\[
   J^*_{\ell, L}(p) = \sup_{a \in X} \sup_{y \in L_X^{\beta}} \sup_{u \in
    A^{\alpha}_{X^*}} \left\{ (a,p(0)) + \int_0^T [\braket{y}{\dot{p}}
     - L(t, u+y,\dot{u})]dt -\ell (u(0)+a,u(T)) \right\}.
\]
Set $
     F \stackrel{\mathrm{def}}{=} \left\{ u \in A^{\alpha}_{X^*}\ ;  \ u
     \in L^{\beta}_X \right\} \subseteq A^{\alpha}_{X^*}.$
Then
\[
     J^*_{\ell,L} (p) \geq \sup_{a \in X} \sup_{y \in L_X^{\beta}} \sup_{u
     \in F} \left\{ (a,p(0)) + \int_0^T [-L(t, u+y,\dot{u}) +
\braket{y}{\dot{p}}]dt
         -\ell (u(0)+a,u(T)) \right\}
\]
Make a substitution $u+y=y' \in L_X^\beta$ to  obtain
\[
     J^*_{\ell,L} (p) \geq \sup_{a \in X} \sup_{y' \in L_X^{\beta}}
     \sup_{u \in F} \left\{ (a,p(0)) - \ell (a+u(0),u(T)) +
     \int_0^T[ \braket{y'-u}{\dot{p}} - L(t, y',\dot{u})] dt  \right\}.
\]
Set now $
S=\{ u:[0,T]\rightarrow X;  u\in L_X^\beta, \, \dot{u}\in L_X^\beta,\,
u(0)\in X\}. $
   Since $\beta\geq 2\geq \alpha$ and ${\| \cdot\|}_{X^*}\leq C\|\cdot\|_X$,
we have
$S\subseteq A_{_{X^*}}^{^\alpha} \cap L_X^\beta =F$ and
\[
     J^*_{\ell ,L} (p) \geq \sup_{a \in X} \sup_{y' \in L_X^{\beta}} \sup_{u
     \in S} \left\{ \big( a,p(0)\big) +
     \int_0^T[ \braket{y'}{\dot p}-\braket{u}{\dot p}-L(t, y',\dot{u})]dt -
\ell \big( a+u(0),u(T)\big)  \right\}
\]
substitute $
u(0) + a = a' \in X$ and write
 \[
     J^*_{\ell ,L} (p) \geq \sup_{a' \in X} \sup_{y' \in L_X^{\beta}} \sup_{u
     \in S} \left\{ \big( a' - u(0),p(0)\big) +
     \int_0^T[ \braket{y'}{\dot p}-\braket{u}{\dot p}-L(t, y',\dot{u})]dt -
\ell \big( a',u(T)\big)  \right\}
\]
Since $\dot u \in L_X^\beta$ and $u\in L_X^\beta$, we have $
\int_0^T \braket{u}{\dot{p}}dt=- \int_0^T \braket{\dot u}{p}dt +
\braket{p(T)}{u(T)} - \braket{p(0)}{u(0)}.
$
 which implies
   \[
     J^*_{\ell ,L} (p) \geq \sup_{a' \in H} \sup_{y' \in
L_X^{\beta}} \sup_{u   \in S}
\left\{ \big( a',p(0)\big) +
     \int_0^T
\left\{\braket{y'}{\dot{p}}+\braket{\dot{u}}{p}-L(t, y',\dot{u})\right\}
dt-\braket{u(T)}{p(T)}-\ell \big( a',u(T)\big)  \right\}.
\]
It is now convenient to identify $S=\{ u: [0,T]\rightarrow X;\,  u\in
L_X^\beta, \, \dot{u}\in L_X^\beta, \, u(0)\in X\}$ with $X\times L_X^\beta$
via the correspondence
$
     (c,v) \in X\times L_X^\beta \mapsto  c+\int_t^T v(s)\, ds\in S$ and 
    $ u\in S   \mapsto  \big( u(T),-\dot u(t)\big)\in X\times L_X^\beta.$ 
 We finally obtain
\begin{eqnarray*}
     J^*_{\ell ,L}(p)&\geq& \sup_{a' \in X} \sup_{c\in X}
\left\{ \big( a',p(0)\big) +  \braket{-c}{p(T)}-\ell (a',c)\right\}\\
&&+\sup_{y'\in L_X^\beta}\sup_{v\in L_X^\beta}\left\{
\int_0^T[\braket{y'}{\dot{p}} +\braket{v}{p} - L(t, y',v)]dt\right\}\\
&=& \int_0^T L^*(t, \dot p (t), p(t)) dt +\ell^* (p(0), -p(T))\\
   &=& \int_0^T L(t, -R p (t), -R\dot p(t)) dt +\ell (-Rp(0), -Rp(T))\\
  &=& I_{\ell, L}(-Rp).
\end{eqnarray*}

\noindent {\bf An application to infinite dimensional Hamiltonian systems:} Let now $Y$ be a reflexive Banach space that is densely embedded in a Hilbert space $E$. Then the product $X := Y \times Y$ is clearly a reflexive Banach space that is densely embedded in the Hilbert space $H = E \times E$. Therefore we have an evolution triple $X \subset H \subset X^*$.\\
We shall consider a simple but illustrative example. Let $\phi_1, \phi_2$ be convex lower semi-continuous functions on $E$ whose domain is $Y$ and is coercive on $Y$. 
Define the convex function $\Phi : H \to \R\cup\{+\infty\}$ by $\Phi(x) = \phi(y_1,y_2) := \phi_1(y_1) + \phi_2(y_2)$.\\
Finally, define the linear automorphism $S : X^* \to X^*$ by $Sx^* = E(y_1^*, y_2^*) := (y_2^*, y_1^*)$. Clearly  S is an automorphism whose restriction to $H$ and $X$ are also automorphisms.\\\\
 Consider now the Lagrangians $L : X^*\times X^* \to \R\cup\{+\infty\}$ defined as:
\begin{eqnarray}
L(x,v) = \Phi(x) + (\Phi\vert_X)^*(-Sv)
\end{eqnarray}
Now for the boundary, consider convex, lower semi-continuous functions $\psi_1$, $\psi_2$: $Y^*\to \R\cup\{\infty\}$ assuming that both are coercive on $Y$. To these functions we associate the boundary Lagrangian $\ell: X^*\times X^* \to \R\cup\{+\infty\}$ by:
\begin{eqnarray}
\label{ell}
\ell((a_1,a_2),(b_1, b_2)) = \psi_1(a_1) + (\psi_1\vert_X)^*(-a_2) +\psi_2(b_1) +(\psi_2\vert_X)^*(b_2)
 \end{eqnarray}
 It is then easy to show that  $L$ is S-anti-selfdual on $X^*\times X^*$ since the convex function $\Phi$ is coercive on $X$ and that $\ell$ is $S$-selfdual.
\begin{proposition}
Suppose that $\phi_j(y) \leq C(\Vert y\Vert_Y^\beta +1)$ for $j = 1,2$, that $\psi_1$ is bounded on the bounded sets of $Y$ and consider the Hamiltonian $\H(p,q)=\phi_1(p) - \phi_2(q)$.   Then for any $T>0$, there exists solutions $(\bar p, \bar q) \in A^\alpha_{H, X^*}$  for the following Hamiltonian system:
\begin{eqnarray*}
\dot p (t) &\in& \partial_2\H(p(t),q(t))\\
-\dot q(t)& \in& \partial_1\H(p(t),q(t))\\
-p(0) \in\partial\psi_1(q(0)) &\&&  p(T) \in \partial\psi_2(q(T)).
\end{eqnarray*}
It can be obtained  by minimizing the following functional on the space $A^\alpha_{H, X^*}$
\[I(p,q) =   \int_0^T \phi((p(t),q(t)) + (\phi_{|_X})^*(-\dot{q}(t),-\dot{p}(t))dt + \ell\big((p(0),q(0)), (p(T),q(T))\big)
\]
where $\phi$ is the convex function $\phi (p,q)=\phi_1(p) + \phi_2(q)$ and $\ell$ is as in (\ref{ell}).
\end{proposition}
\noindent{\bf Proof:}
We wish to apply Theorem 2.7  to the S-anti-selfdual Lagrangian pair $(L,\ell)$ defined above, so we must check that $J^\alpha_{\ell ,L}:
H\times L_X^{\beta} \rightarrow {\R \cup \{+\infty\}}$ is
sub-differentiable at
$(0,0)$. 
To do this we use the assumption on $\phi_j$ to obtain the inequality:
\begin{eqnarray*}
J^\alpha_{\ell ,L}(a,v) &=&
 \inf \{ \int_0^T L(t, u+v,\dot{u})dt + \ell
(u(0)+a,u(T))\ ;
\ u \in A^{\alpha}_{X^*} \}\\
&\leq&\int_0^T L(t, v,0)dt + \ell(a,0)\\
&\leq&C(\frac{\Vert a\Vert^2_H}{2} + \int_0^T {\Vert v\Vert^\beta_X dt}+1).
\end{eqnarray*}
Again, since $J^\alpha_{\ell ,L}$ is bounded on bounded sets of $H\times L_X^{\beta}$, we conclude that it is subdifferentiable at $(0,0)$. Thus  there exists $\bar x(\cdot) = \big(\bar p(\cdot), \bar q(\cdot)\big)\in A^{\alpha}_{X^*}$ such that $ I_{L,\ell}(\bar x(\cdot)) = 0$. Therefore,
\begin{eqnarray*}
0 &=& \int_0^T L(t, \bar x,\dot{\bar x})dt + \ell
(\bar x(0), \bar x(T))\\
& = &\int_0^T \phi(\bar x) + (\phi\vert_X)^*(-S\dot{\bar x})dt + \ell
(\bar x(0), \bar x(T))\\
&\geq& -\int_0^T {\langle\bar x,S\dot{\bar x}\rangle dt} + \ell
(\bar x(0), \bar x(T))\\
& = &-\int_0^T {\frac{d\langle \bar p, \bar q\rangle}{dt} dt} + \psi_1(\bar p(0)) + (\psi_1\vert_X)^*(-\bar q(0)) + \psi_2(\bar p(T)) + (\psi_2\vert_X)^*(\bar q(T)) \\
&=& \braket{\bar p(0)}{\bar q(0)} - \braket{\bar p(T)}{\bar q(T)} +\psi_1(\bar p(0)) + (\psi_1\vert_X)^*(-\bar q(0)) + \psi_2(\bar p(T)) + (\psi_2\vert_X)^*(\bar q(T))  \\ &\geq& 0.
\end{eqnarray*}
 Therefore every inequality in this chain is actually an equality. We conclude   that
$-S\dot{\bar x}(t) \in \partial\Phi(\bar x(t))$ for almost all $t\in[0,T]$
 and that
\[\braket{\bar p(0)}{\bar q(0)}+\psi_1(\bar p(0)) + (\psi_1\vert_X)^*(-\bar q(0)) = - \braket{\bar p(T)}{\bar q(T)} + \psi_2(\bar p(T)) + (\psi_2\vert_X)^*(\bar q(T))=0\]
 
By the definition of $S$ and $\Phi$ and Fenchel inequality, this is precisely a solution of the equation above. 
 
 \section{Two-parameter gradient flows}

Behind  the results of the previous section is the fact that an R-antiselfdual  Lagrangian on a Hilbert space $H$ lifts to an R-antiselfdual  Lagrangian on path space. So far, we only needed anti-selfduality  on the elements of $A^2_H\times \{0\}$. However, we have the following stronger stability result announced in \cite{G01} and proved in \cite{G3}. For clarity we shall restrict ourselves to ASD-Lagrangians (i.e., $R(x)=x$).

 \begin{lemma}
\label{iterate ASD lagrangian}
Let $H$ be a Hilbert space and let $L : [0,T]\times H\times H\to \R$ be an anti-selfdual Lagrangian  such that for every $p \in H$ and $t\in [0,T]$ the map
\begin{equation}
\label{good.riddance}
x\mapsto L(t,x,p)\quad \mbox{is bounded on the bounded sets of $H$.}
\end{equation}
 Then for every $x_0\in H$, the Lagrangian defined on $ L^2_H([0,T]) \times L^2_H([0,T])$   by
\[{\cal L}(x,p) := 
\begin{cases}
{\int_0^T L(t,x(t),\frac{dx}{dt}(t) + p(t))dt + \frac{1}{2} \|x(0)\|^2_H + 2\braket{x_0}{x(0)} + \|x_0\|^2_H +  \frac{1}{2}\|x(T)\|_H^2\mbox{ if $x(\cdot)\in A^2_{H'}$}}\cr{\infty\quad \quad\quad\quad\quad\quad\quad\quad\quad\quad\quad\quad\quad\quad\quad\quad\quad\quad\quad\quad\quad\quad\quad\quad\quad\quad\quad\quad\quad\quad\ \ \mbox{ otherwise}}\cr
\end{cases}
\]
 is also an ASD Lagrangian on $ L^2_H([0,T]) \times L^2_H([0,T])$.
\end{lemma}
\noindent{\bf Proof:} Note that this also follows from a more general result established in \cite{GT2}. Indeed, since $L(t,x,p)$ is an anti-self-dual Lagrangian on $H$, the map
\[(x(\cdot),p(\cdot))\mapsto\int_0^T L(t, x(t),p(t))dt\]
is an ASD Lagrangian on the path space $L^2_H([0,T]) \times L^2_H([0,T])$ (See \cite{G3}). Now, using the terminology of \cite{GT2},  the map $x \mapsto \frac{dx}{dt}$ (with domain $A^2_H([0,T])$) is  skew-adjoint modulo the boundary operator $x\to (x(0), x(T))$  on the Hilbert space $L^2_H([0,T])$. Therefore
 ${\cal L}$ is also an ASD Lagrangian. \\

 Setting $H' :=L^2_H([0,T])$ as a state space, and since  ${\cal L}(\cdot,\cdot) : H'\times H'\to\R$  is now an  anti-selfdual Lagrangian on $H'$, we can then lift it to a new path space  $L^2_{H'}([0,S])$ and obtain a new action functional 
 \[
{\cal I}(x(\cdot)) := \int_0^S{\cal L}(x(s), \frac{dx}{ds}(s))ds + \ell'(x(0), x(S))
 \]
that we can minimize on  $A^2_{H'}([0,S])$. Here is the main result of this section. We recall from \cite{GT2} that the partial domain ${\rm Dom}_1(\partial L)$ of a Lagrangian $L$ is defined as:
\[
{\rm Dom}_1(\partial L)=\{x\in H; \hbox{\rm There exists $p\in H\ \mbox{such that }\ -(p,x)\in \partial L(x, p)$}\}.
\]
\begin{theorem}
\label{main theorem}
Let $H$ be a Hilbert space and $L :  H\times H\to \R$ be an ASD Lagrangian that is uniformly convex in the first variable. If $x_0\in {\rm Dom}_1(\partial L)$, then there exists $\hat x(\cdot,\cdot)\in A^2([0,S]; L^2_H([0,T])$ such that $\hat x(s,\cdot)\in A^2_H([0,T])$ for almost all $s\in [0,S]$ and
\begin{eqnarray}
\label{minimum is zero}
0&=&\int_0^S \int_0^T L(\hat x(s,t),\frac{\partial\hat x}{\partial t}(s,t) + \frac{\partial \hat x}{\partial s}(s,t))dtds\nonumber\\
&& + \int_0^S \big(\frac{1}{2} \|\hat x(s,0)\|^2_H -2\braket{\hat x(s,0)}{x_0} + \|x_0\|_H^2 +  \frac{1}{2}\|\hat x(s,T)\|_H^2\big)ds\nonumber\\
&& + \int_0^T \big(\frac{1}{2} \|\hat x(0,t)\|^2_H -2\braket{\hat x(0,t)}{x_0} + \|x_0\|_H^2 +  \frac{1}{2}\|\hat x(S,t)\|_H^2\big)dt.
\end{eqnarray}
Furthermore, for almost all $(s,t)\in[0,S]\times[0,T]$, we have
\begin{eqnarray}
\label{main equation}
-\frac{\partial \hat x}{\partial t}(s,t) - \frac{\partial \hat x}{\partial s}(s,t)&\in& \partial_1 L(\hat x(s, t), \frac{\partial \hat x}{\partial t}(s,t) + \frac{\partial \hat x}{\partial s}(s,t))\\
 - \hat x(s,t) &\in&   \partial_2L(\hat x(s, t), \frac{\partial \hat x}{\partial t}(s,t) + \frac{\partial \hat x}{\partial s}(s,t))\\
\hat x(0,t) &=& x_0\ a.e.\ t\in [0,T]\\
\hat x(s,0) &=& x_0\ a.e.\ s\in [0,S].
\end{eqnarray}
\end{theorem}
We first note that if $L$ satisfies the boundedness condition (\ref{good.riddance}) then the conclusions of the theorem are easy to establish as shown in the following Lemma. The main difficulty of the proof is to get rid of this 
condition.

 \begin{lemma}
\label{existence assuming boundedness}
Let $L : [0,T]\times H\times H\to \R$ be an ASD Lagrangian  on a Hilbert space $H$ such that $L(t, \cdot, \cdot)$  is uniformly convex in the first variable for each $t\in [0, T]$ while verifying condition  (\ref{good.riddance}). If $x_0\in {\rm Dom}_1(\partial L)$ then there exists $\hat x(\cdot,\cdot)\in A^2([0,S]; L^2_H([0,T])$ such that $\hat x(s,\cdot)\in A^2_H([0,T])$ for almost all $s\in [0,S]$ and satisfying properties  (12)-(15) above. 
 \end{lemma}
\noindent{\bf Proof:}
According to Lemma \ref{iterate ASD lagrangian},  ${\cal L}$ is a uniformly convex ASD Lagrangian on $H' := L^2_H([0,T])$. Since $0\in {\rm Dom}_1(\partial L)$ we have that $0\in {\rm Dom}_1(\partial {\cal L})$. Therefore by Theorem 4.1 of \cite{GT2},  we can find an
$\hat x(\cdot)\in A^2_{H'}([0,S]) = A^2([0,S]; L^2_H([0,T]))$ such that
\[0=\int_0^T {\cal L}(\hat x(t), \dot{\hat x}(t))dt + \frac{1}{2}\|\hat x(0)\|_{H'}^2 +\frac{1}{2}\|\hat x(T)\|_{H'}^2.\]
From the definintion of  ${\cal L}$, we get that $\hat x(s,\cdot)\in A^2_H([0,T])$ for almost all $s\in [0,S]$ while satisfying (\ref{minimum is zero}).  We therefore get the following chain of inequalities:
\begin{eqnarray*}
0&=&\int_0^S \int_0^T L(t,\hat x(s,t),\frac{\partial\hat x}{\partial t}(s,t) + \frac{\partial\hat x}{\partial s}(s,t))dtds\\
&& + \int_0^S \big(\frac{1}{2} \|\hat x(s,0)\|^2_H -2\braket{\hat x(s,0)}{x_0} + \|x_0\|_H^2 +  \frac{1}{2}\|\hat x(s,T)\|_H^2\big)ds\\
&& + \int_0^T\big( \frac{1}{2} \|\hat x(0,t)\|^2_H -2\braket{\hat x(0,t)}{x_0} + \|x_0\|_H^2 +  \frac{1}{2}\|\hat x(S,t)\|_H^2\big)dt\\
&\geq&\int_0^S\int_0^T-\braket{x(s,t)}{\frac{\partial\hat x}{\partial t}(s,t) + \frac{\partial\hat x}{\partial s}(s,t)}dtds\\&&+ \int_0^S\big( \frac{1}{2} \|\hat x(s,0)\|^2_H -2\braket{\hat x(s,0)}{x_0} + \|x_0\|_H^2 +  \frac{1}{2}\|\hat x(s,T)\|_H^2\big)ds\\
&& + \int_0^T\big( \frac{1}{2} \|\hat x(0,t)\|^2_H -2\braket{\hat x(0,t)}{x_0} + \|x_0\|_H^2 +  \frac{1}{2}\|\hat x(S,t)\|_H^2\big)dt\\
&\geq& \int_0^T\|\hat x(0,t)-x_0\|_H^2 dt + \int_0^S\|\hat x(s,0)-x_0\|_H^2 ds\geq 0.
\end{eqnarray*}
This means that for almost all  $(s,t)\in[0,S]\times[0,T]$
 \[-\frac{\partial\hat x}{\partial t}(s,t) - \frac{\partial\hat x}{\partial s}(s,t)\in\partial_1 L(t,\hat x(t), \frac{\partial\hat x}{\partial t}(s,t) + \frac{\partial\hat x}{\partial s}(s,t))\]
\[\hat x(s,t) \in\partial_2L(t,\hat x(t), \frac{\partial\hat x}{\partial t}(s,t) + \frac{\partial\hat x}{\partial s}(s,t))\]
\[\hat x(0,t) = x_0\ {\rm a.e.} \ t\in [0,T]\]
\[\hat x(s,0) = x_0\ {\rm a.e.} \ s\in [0,S]\]\\
In the next proposition we do away with the assumption of boundeness of the ASD Lagrangian $L$ that was used in Lemma \ref{existence assuming boundedness}. The argument we use is similar to that in \cite{GT2}.  We first $\lambda$-regularize the Lagrangian $L$ then derive some uniform bounds to ensure convergence in the proper topology when $\lambda$ goes to 0. To do this we need to first state some precise estimates on approximate solutions obtained using inf-convolution. Recall first from \cite{G2} that the Lagrangian
\[
L_\lambda(x,p) := \inf\limits_{z\in H}\lbrace L(z,p) + \frac{1}{2\lambda}\|x-z\|^2_H\rbrace + \frac{\lambda}{2}\|p\|^2_H
\]
is anti-selfdual  for each $\lambda>0$. 
 \begin{lemma} For a given convex functional $L:H\times H\to \R\cup\{+\infty\}$ and  $\la >0$, denote for each $(p,x)\in H\times H$, by $J_\la (x,p)$ the minimizer of the following optimization problem:
\[
 \inf_z\left\{ L(z,p)+\frac{\| x-z\|_H^2}{2\la}\right\}.
 \]
Then for each $(x,p)\in H\times H$, we have
 \begin{equation}
\label{Jlambda}
\partial_1L_\la (x,p)=\frac{x-J_\la (x,p)}{\la}\in \partial_1L\big( J_\la (x,p),p\big).
\end{equation}
 \end{lemma}
\noindent{\bf Proof:} This is left to the reader.

\begin{lemma}
\label{ASD bound}
Assume  $L:H\times H\to\overline{\R}$  is an anti-selfdual Lagrangian and let $L_\lambda$ be its $\lambda$-regularization, then the following hold:
\begin{enumerate}
\item If $ -(y,x)=\partial L_\la (x,y)$, then necessarily
$-\big( y,J_\la (x,y)\big)\in \partial L\big( J_\la (x,y),y\big).$
\item If $x_0\in {\rm Dom}_{1}(\partial L)$, then  $\| y_\la\|_H\le \|\hat p\|_H$ whenever  $y_\la$ solves 
$
-(y_\la ,x_0)=\partial L_\la (x_0,y_\la )
$ and $\hat p$ solves $-(\hat p ,x_0)=\partial L (x_0,\hat p )$.
\end{enumerate}
\end{lemma}
\noindent {\bf Proof:}  (1) \quad If $ -(y,x)=\partial L_\la (x,y)$ then  $ L_\la (x,y)+L_\la ^*(-y,-x)=-2\langle x,y)$ and since 
$L$ is an ASD Lagrangian, we have $L_\la (x,y)+L_\la (x,y)=-2\langle x,y)$, hence 
 \begin{eqnarray*}
-2\braket{x}{y} &=& L_\la (x,y)+L_\la (x,y)\\
 &=&2\left( L\big( J_\la (x,y),y\big) 
+\frac{\| x-J_\la (x,y)\|_H^2}{2\la}+\frac{\la\| y\|_H^2}{2}\right)\\
&=& L^*\big( -y,-J_\la (x,y)\big) +L\big( J_\la (x,y),y\big) 
   +2\left(\frac{\|-x+J_\la (x,y)\|_H^2}{2\la}+\frac{\la\| y\|_H^2}{2}\right)\\
&\ge& -2\langle y,J_\la (x,y)\rangle +2 \langle -x+J_\la (x,y),y\rangle\\
&=& -2\langle x,y).
\end{eqnarray*}

The second last inequality is deduced by applying Fenchel's inequality to the first two terms and the last two terms. The above chain of inequality shows that all inequalities are equalities. This implies, again by Fenchel's inequality that $-\big( y,J_\la (x,y)\big)\in \partial L\big( J_\la (x,y),y\big)$.\\
  
(2) If $-(y_\la ,x_0)=\partial L_\la (x_0,y_\la )$, we get from the previous lemma that 
\[
-y_\la =\frac{x_0-J_\la (x_0,y_\la )}{\la}\in \partial_1L\big( J_\la (x_0,y_\la ),y_\la\big),
\]
 and by the first part of this lemma, that 
 \[
 -\big( y_\la ,J_\la (x_0,y_\la )\big)\in \partial L\big( J_\la (x_0,y_\la ),y_\la\big).
 \] 
 Now since $x_0\in {\rm Dom}_{1}(\partial L)$, there exists $\hat p$ such that $(-\hat p, -x_0)\in  \partial L(x_0,\hat p)$. 
Setting $ v_\la = J_\la (x_0,y_\la )$, and 
since  $-\big( y_\la ,v_\la )\big)\in \partial L\big( v_\la, y_\la\big)$, we get from monotonicity and by the fact that $y_\la=\frac{v_\la - x_0}{\la}$,
\begin{eqnarray*}
0 &\le&\braket{(x_0,\hat p)-(v_\la ,y_\la )}{\big( \partial_1L(x_0,\hat p),\partial_2L(x_0,\hat p)\big) -(-y_\la ,-v_\la )}\\
&=&\braket{(x_0,\hat p)-(v_\la ,y_\la )}{(-\hat p,-x_0)-(\frac{x_0-v_\la}{\la},-v_\la )}\\
  &=&-\frac{\|x_0 - v_\la\|_H^2}{\la}+\braket{v_\la - x_0}{\hat p} 
    +\braket{\hat p}{v_\la-x_0} -\braket{y_\la}{v_\la-x_0}\\
  &=&-2\frac{\|x_0- v_\la\|_H^2}{\la}+2\braket{v_\la-x_0}{\hat p}
\end{eqnarray*} 
which yields that $\frac{\| x_0-v_\la\|_H}{\la}\le \|\hat p\|_H$ and finally the desired bound
$\| y_\la\|\le\|\hat p\|$ for all $\la >0$. \\
\begin{lemma}
\label{bound by p_0}
Let $L : H\times H$ be an anti-self-dual Lagrangian that is uniformly convex in the first variable. If $x_0\in {\rm Dom}_1(\partial L)$ and if $\hat x(\cdot)\in A^2_H([0,T]$ satisfies
\[\int_0^T L(\hat x(t),\dot{\hat x}(t))dt + \frac{1}{2}\|\hat x(0)\|_H^2 + \|x_0\|_H^2 +\braket{\hat x(0)}{x_0} + \frac{1}{2}\|\hat x(T)\|_H^2 = 0, \]
then we have the estimate
\begin{equation}
\int_0^T\|\dot{\hat x}(t)\|_H^2dt \leq T\|p_0\|_H^2,
\end{equation}
  where $p_0$ is the point that satisfies $-(p_0,x_0)\in \partial L(x_0,p_0)$
\end{lemma}
\noindent{\bf Proof:}
By the uniquenss of the minimizer established in \cite{GT2}, $\hat x(\cdot)$ is the weak  limit in $A^2_H([0,T])$ of $\{x_\la(\cdot)\in C^{1,1}([0,T])\}$ where $(-\dot x_\la(t),-x_\la(t))\in\partial L_\la(x_\la(t),\dot x_\la(t))$, $x_\la(0) =x_0$.
 
Standard arguments using monotonicity shows that $\|\dot x_\la(t)\|_H\leq \|\dot x_\la(0)\|_H$ for all $t\in[0,T]$. Since $(-\dot x_\la(0), -x_0)\in\partial L_\la(x_0, \dot x_\la(0))$, Lemma \ref{ASD bound} shows that $\|x_\la(0)\|_H\leq \|p_0\|_H$ for all $\la >0$. Therefore, letting $\la\to0$ and by weak lower semi-continuity of the norm we get that $\int_0^T\|\dot{\hat x}(t)\|^2_H dt\leq T\|p_0\|^2_H$.\\

\noindent{\bf Proof of Theorem 3.2:}
Apply Lemma \ref{existence assuming boundedness} to $L_\lambda$ we obtain an $\hat x_\lambda(\cdot)\in A^2([0,S];L^2_H([0,T])$ satisfying for all $(s,t)\in[0,S]\times[0,T]$
\begin{eqnarray}
\label{lambda equation}
 -\frac{d\hat x_\lambda}{dt}(s,t) - \frac{d\hat x_\lambda}{ds}(s,t)\in \partial_1 L_\lambda(\hat x_\lambda(t), \frac{d\hat x_\lambda}{dt}(s,t) + \frac{d\hat x_\lambda}{ds}(s,t))\nonumber\\
-\hat x_\lambda(s,t)\in\partial_2L_\lambda(\hat x_\lambda(t), \frac{d\hat x_\lambda}{dt}(s,t) + \frac{d\hat x_\lambda}{ds}(s,t))
\end{eqnarray}
\[\hat x_\lambda(0,t) = x_0\ \forall\ t\in [0,T]\]
\[\hat x_\lambda(s,0) = x_0\ \forall\ s\in [0,S]\]
and
\begin{eqnarray}
\label{approximate minimizer}
0&&=\int_0^S \int_0^T L_\lambda(\hat x_\lambda(s,t),\frac{d\hat x_\lambda}{dt}(s,t) + \frac{d\hat x_\lambda}{ds}(s,t))dtds\nonumber\\
&& + \int_0^S \big(\frac{1}{2} \|\hat x_\lambda(s,0)\|^2_H -2\braket{x_\lambda(s,0)}{x_0} + \|x_0\|_H^2 +  \frac{1}{2}\|\hat x_\lambda(s,T)\|_H^2\big)ds\nonumber\\
&& + \int_0^T \big(\frac{1}{2} \|\hat x_\lambda(0,t)\|^2_H -2\braket{x_\lambda(0,t)}{x_0} + \|x_0\|_H^2 +  \frac{1}{2}\|\hat x_\lambda(S,t)\|_H^2\big)dt.
\end{eqnarray}
Now consider the ASD Lagrangian ${\cal L}_\lambda$ on $L_H^2([0,T])$ defined by:
{\small
\begin{eqnarray*}
{\cal L}_\lambda(x,p) := \begin{cases}
{\int_0^T L_\lambda(x(t),\frac{dx}{dt}(t) + p(t))dt + \frac{1}{2} \|x(0)\|^2_H +  \frac{1}{2}\|x(T)\|_H^2 - 2\braket{x_0}{x(0)}+ \|x_0\|^2_H \mbox{ if $x \in A^2_{H}([0,T])$}}\cr{\infty\mbox{\quad\quad\quad\quad\quad\quad\quad\quad\quad\quad\quad\quad\quad\quad\quad\quad\quad\quad\quad\quad\quad\quad\quad\quad\quad\quad\quad\quad\quad\quad\quad\quad\quad\quad else}}\cr
\end{cases}
\end{eqnarray*}
}
 
Let $\hat{\cal X}_\lambda(\cdot) : [0,S]\to L^2_H([0,T])$ be the map $s\mapsto \hat x_\lambda(s,\cdot)\in L^2_H([0,T])$ and denote by ${\cal X}_0\in L^2_H([0,T])$ the constant map $t\mapsto x_0$. Then by (\ref {approximate minimizer}) $\hat{\cal X}_\lambda(\cdot)$ is the arc in $A^2([0,S]; L^2_H([0,T]))$ satisfying
\[0 = \int_0^S {\cal L}_\lambda(\hat{\cal X}_\lambda(s), \frac{d\hat{\cal X}_\lambda}{ds}(s))ds + \frac{\|\hat{\cal X}_\lambda(0)\|^2_{L^2_H([0,T])}}{2}+\frac{\|\hat{\cal X}_\lambda(S)\|^2_{L^2_H([0,T])}}{2} -2\braket{{\cal X}_0}{\hat{\cal X}_\lambda(S)}_{L^2_H([0,T])}  + \|{\cal X}_0\|_{L^2_H([0,T])}\]
with ${\cal X}_0 \in {\rm Dom}_1(\partial{\cal L}_\lambda)$. Apply Lemma \ref{bound by p_0} to the ASD Lagrangian ${\cal L}_\lambda$ and the Hilbert space $L^2_H([0,T])$ we get that 
\[\int_0^S\int_0^T \|\frac{d\hat x_\la(s,t)}{ds}\|_H^2dtds \leq S\int_0^T\|{\cal P}_\lambda(t)\|^2_Hdt\]
where ${\cal P}_\lambda\in L^2_H([0,T])$ is any arc that satisfies $(-{\cal P}_\lambda,-{\cal X}_0)\in\partial{\cal L}_\lambda({\cal X}_0, {\cal P}_\lambda)$. Observe that if the point $p_\lambda \in H$ satifies the equation $-(p_\lambda, x_0)\in\partial L_\lambda(x_0,p_\lambda)$, then we can just take ${\cal P}_\lambda$ to be the constant arc $t\mapsto p_\lambda$. Combining this fact with Lemma \ref{ASD bound}, we obtain that for all $s\in [0,S]$ and all $\lambda >0$, 
\[\int_0^S\int_0^T \|\frac{d\hat x_\la(s,t)}{ds}\|_H^2dtds \leq ST\|p_0\|_H^2\]
In deriving the above estimates, we interpreted $\hat x_\lambda(s,t)$ as a map $\hat{\cal X}_\lambda(\cdot) :[0,S]\to L^2_H([0,T])$. However, we can also view it as a map from $[0,T]\to L^2_H([0,S])$ and run the above argument in this new setting. By doing this we obtain that for all $\lambda >0$:
\begin{eqnarray}
\label{uniform estimate on derivatives}
\int_0^S\int_0^T \|\frac{d\hat x_\la(s,t)}{ds}\|_H^2dtds +\int_0^S\int_0^T \|\frac{d\hat x_\la(s,t)}{dt}\|_H^2dtds \leq 2TS\|p_0\|_H^2.
\end{eqnarray}
 Now for any $\big(v_1(s,t),\ v_2(s,t)\big)$ satisfying equation (\ref{lambda equation}) we can use monotonicity to derive the bound:
\begin{eqnarray*}
&&\frac{d}{dt}\|v_1(s,t) - v_2(s,t)\|_H^2 + \frac{d}{ds}\|v_1(s,t) - v_2(s,t)\|_H^2\leq 0.
\end{eqnarray*}
So we obtain
\begin{eqnarray*}
&&\int_0^S\|v_1(s,t) - v_2(s,t)\|_H^2ds + \int_0^T\|v_1(s,t) - v_2(s,t)\|_H^2dt\\
&&\leq \int_0^S\|v_1(s,0) - v_2(s,0)\|_H^2ds+\int_0^T\|v_1(0,t) - v_2(0,t)\|_H^2dt.
\end{eqnarray*}
Now pick $v_1(s,t) = \hat x_\lambda(s,t)$ and $v_2(s,t) = \hat x_\lambda(s+h,t)$ we get that
\begin{eqnarray}
\label{boundedness by initial data}
\int_0^S\|\frac{d}{ds}\hat x_\lambda(s,t)\|_H^2ds +\int_0^T\|\frac{d}{ds}\hat x_\lambda(s,t)\|_H^2dt\leq
\int_0^S\|\frac{d}{ds}\hat x_\lambda(s,0)\|_H^2 ds+\int_0^T\|\frac{d}{ds}\hat x_\lambda(0,t)\|_H^2dt.
\end{eqnarray}
Setting $s = 0$ in equation (\ref {lambda equation}) we get that for all $t\in [0,T]$
\[-\big(\frac{d\hat x_\lambda}{dt}(0,t) + \frac{d\hat x_\lambda}{ds}(0,t), x_0\big)\in \big(\partial_1 L_\lambda(x_0, \frac{d\hat x_\lambda}{dt}(0,t) + \frac{d\hat x_\lambda}{ds}(0,t)), \partial_2L_\lambda(x_0, \frac{d\hat x_\lambda}{dt}(0,t) + \frac{d\hat x_\lambda}{ds}(0,t))\big)\]

Therefore by Lemma \ref{ASD bound} we have that for all $t\in [0,T]$ and $\lambda >0$, 
\[
\|\frac{d\hat x_\lambda}{dt}(0,t) + \frac{d\hat x_\lambda}{ds}(0,t)\|_H \leq \|p_0\|_H.
\]
 Observe that if we take $v_2(s,t) = \hat x_\lambda(s,t+h)$ we can use the same argument as above to get that for all $s\in[0,S]$, 
 \[
 \|\frac{d\hat x_\lambda}{dt}(s,0) + \frac{d\hat x_\lambda}{ds}(s,0)\|_H \leq \|p_0\|_H.
 \]
   Therefore, for all $s\in [0,S]$, $t\in [0,T]$, and $\lambda >0$:
\begin{eqnarray}
\label{bound on initial derivatives}
\|\frac{d\hat x_\lambda}{dt}(0,t) + \frac{d\hat x_\lambda}{ds}(0,t)\|_H + \|\frac{d\hat x_\lambda}{dt}(s,0) + \frac{d\hat x_\lambda}{ds}(s,0)\|_H\leq 2\|p_0\|_H.
\end{eqnarray}
Combining (\ref{bound on initial derivatives}), (\ref {boundedness by initial data}), and (\ref{uniform estimate on derivatives}) we get that
\begin{eqnarray}
\label{L2 bound on derivatives}
\int_0^S\int_0^T\|\frac{d\hat x_\lambda}{ds}(s,t)\|_H^2 +\|\frac{d\hat x_\lambda}{dt}(s,t)\|_H^2dtds \leq C,
\end{eqnarray}
for some constant independent of $\lambda$. If we denote by $J_\lambda(x,p)$ the point that satisfies $L_\lambda(x,p) = L(J_\lambda(x,p),p) + \frac{\lambda}{2}\|p\|_H^2$ and $v_\lambda(s,t)$ to be $J_\lambda\big(\hat x_\lambda(s,t), \frac{d\hat x_\lambda}{dt}(s,t) + \frac{d\hat x_\lambda}{ds}(s,t)\big)$, then we can deduce from equation (\ref{lambda equation}) that
\[-\frac{d\hat x_\lambda}{dt}(s,t) - \frac{d\hat x_\lambda}{ds}(s,t) = \frac{\hat x_\lambda(s,t) - v_\lambda(s,t)}{\lambda}.
\]
The estimate given by equation (\ref{L2 bound on derivatives}) then implies
\[\lim\limits_{\lambda\to 0}\int_0^T\int_0^S \|\hat x_\lambda(s,t) - v_\lambda(s,t)\|^2_Hdsdt = 0\]
Therefore, combining this with (\ref{L2 bound on derivatives}) we obtain the following convergence result:
\begin{eqnarray}
\label{convergence}
\hat x_\lambda (\cdot,\cdot) \rightharpoonup \hat x(\cdot,\cdot)\ \mbox{in $A^2([0,S]; L^2_H([0,T])$}\\
\hat x_\lambda (\cdot,\cdot) \rightharpoonup \hat x(\cdot,\cdot)\ \mbox{in $A^2([0,T]; L^2_H([0,S])$}\\
v_\lambda(\cdot,\cdot) \rightharpoonup \hat x(\cdot,\cdot)\ \mbox{in $L^2_H([0,S]\times[0,T])$.}
\end{eqnarray}
Write (\ref{approximate minimizer}) in the form
\begin{eqnarray*}
0&&=\int_0^S \int_0^T L(v_\lambda(s,t),\frac{d\hat x_\lambda}{dt}(s,t) + \frac{d\hat x_\lambda}{ds}(s,t)) + \frac{\lambda}{2}\|\frac{d\hat x_\lambda}{dt}(s,t) + \frac{d\hat x_\lambda}{ds}(s,t)\|_H^2dtds\\
&& + \int_0^S \frac{1}{2} \|\hat x_\lambda(s,0)\|^2_H -2\braket{x_\lambda(s,0)}{x_0} + \|x_0\|_H^2 +  \frac{1}{2}\|\hat x_\lambda(s,T)\|_H^2ds\\
&& + \int_0^T \frac{1}{2} \|\hat x_\lambda(0,t)\|^2_H -2\braket{x_\lambda(0,t)}{x_0} + \|x_0\|_H^2 +  \frac{1}{2}\|\hat x_\lambda(S,t)\|_H^2dt
\end{eqnarray*}
and take $\lambda \to 0$ using the convergence results in (\ref{convergence}) in conjunction with lower-semi-continuity we get
\begin{eqnarray*}
0&&\geq\int_0^S \int_0^T L(\hat x(s,t),\frac{\partial\hat x}{\partial t}(s,t) + \frac{\partial\hat x}{\partial s}(s,t)) dtds\\
&& + \int_0^S \frac{1}{2} \|\hat x(s,0)\|^2_H -2\braket{\hat x(s,0)}{x_0} + \|x_0\|_H^2 +  \frac{1}{2}\|\hat x(s,T)\|_H^2ds\\
&& + \int_0^T \frac{1}{2} \|\hat x(0,t)\|^2_H -2\braket{\hat x(0,t)}{x_0} + \|x_0\|_H^2 +  \frac{1}{2}\|\hat x(S,t)\|_H^2dt\geq 0
\end{eqnarray*}

Standard arguments then give the desired result.\\\\
Clearly, this argument can be extended to obtain N-parameter gradient flow. We state the result without proof.
\begin{corollary}
Let $L(\cdot,\cdot) : H\times H\to\R\cup\{+\infty\}$ be an ASD Lagrangian that is uniformly convex in the first variable and let $u_0\in {\rm Dom}_1(\partial L)$. Then for all $T_1 \geq T_2..\geq T_N>0$, there exists $u\in  L^2_H(\prod\limits_{j = 0}^N[0, T_j])$ such that $\frac{\partial u}{\partial t_j} \in L^2_H(\prod\limits_{j = 0}^N[0, T_j])$ for all $j = 1,..., N$ and which satisfies the differential equation
\[-\sum\limits_{j = 1}^N \frac{\partial u}{\partial t_j}(t_1,...,t_N)\in\partial_1L(u(t_1,..,t_N),\sum\limits_{j = 1}^N \frac{\partial u}{\partial t_j}(t_1,...,t_N))\] 
with boundary data
$u(t_1,...,t_N)= u_0$ if  one of the $t_j = 0$
\end{corollary} 
We conclude this paper with some remarks.
\begin{remark}\rm 
Let $u:[0,T] \to H$ be the 1-parameter gradient flow associated to an ASD Lagrangian  $L$ (See \cite{GT2}). Namely,
\[-\frac{du}{dt}(t) \in \partial_1L(u(t),\frac{du}{dt}(t))\]
\[u(0)= u_0\]
If we make the change of variables $v(s',t') = u(s'+t')$, then $v(\cdot,\cdot)$ obviously solves (\ref{main equation}), with however the boundary condition $v(s',t') = u_0$ on the hyperplane $s' = -t'$.   In comparison, Theorem \ref{main theorem} above yields a solution $u(\cdot, \cdot)$ for (\ref{main equation}) with a  boundary condition that is prescribed  on two hyperplanes, namely  $u(0,t)=u(s,0) = u_0$ for all $(s,t)\in[0,S]\times[0,T]$. 
\end{remark}
\begin{remark}\rm 
Suppose now $u(\cdot,\cdot):[0,\infty)\times[0,\infty)\to H$ solve (\ref {main equation}) with initial boundary condition $u(0,t) = u(s,0) = u_0$ for all $(s,t) \in[0,\infty)\times[0,\infty)$, and consider the change of variable 
\[
v(s',t') = u(s',(1-C)s' + Ct')
\]
 for some $C >0$. Then $v(s',t')$ again solves (\ref{main equation}) on the domain 
 \[
D= \{(s',t')\in \R\times \R; \, s'\geq 0, t'\geq (1-\frac{1}{C})s'\}.
 \]
  The boundary condition for $v(s',t')$ is 
  \[
\hbox{$  v(0,t') = v(s', \frac{1-C}{-C}s') = u_0$ for all $t'\geq 0$ and $s'\geq 0$.}
  \]
This is  essentially a two-parameter ASD flow on the wedge $D$.
\end{remark}
\begin{remark}\rm 
Let now $u(\cdot,\cdot,\cdot) :[0,\infty)\times[0,\infty)\times[0,\infty)\to H$ be a solution for the three-parameter ASD flow.
\[-\frac{\partial u}{\partial r}-\frac{\partial u}{\partial s}-\frac{\partial u}{\partial t}(r,s,t) \in\partial_1L(u(r,s,t),\frac{\partial u}{\partial r}+\frac{\partial u}{\partial s}+\frac{\partial u}{\partial t}(r,s,t))\]
\[u(0,s,t) = u(r,0,t) = u(r,s,0) = u_0\]
With  the change of variable $v(r',s',t') = u(\frac{s'+r'}{2}, \frac{t'+r'}{2}, \frac{s'+t'}{2})$, $v(r',s',t')$ again solves the differential equation
\[-\frac{\partial v}{\partial r'}-\frac{\partial v}{\partial s'}-\frac{\partial v}{\partial t'} \in\partial_1L(u,\frac{\partial v}{\partial r'}+\frac{\partial v}{\partial s'}+\frac{\partial v}{\partial t'})\]
on the domain 
\[
D=\{(r',s',t')\mid s'\geq -r',\ r'\geq -t',\ s'\geq -t'\}
\]
 with boundary conditions 
 \[
\hbox{$ v(r',s',t') = u_0$ if $s'= -r'\ or\ r' = -t'\ or\ s' = -t'.$}
 \]
 Looking now at $(r',s')$ as "state" variables and $t'$ as the time variable, we see that at any given time $t'$, $v(r',s',t')$ solves the equation on $\{(r',s')\mid s'\geq -r',\ r'\geq -t',\ s'\geq -t'\}$ with $v = u_0$ on the boundary of this domain. This essentially describes a simple PDE with a time evolving boundary.
\end{remark}

\end{document}